\documentclass[11pt, a4paper]{article}
\usepackage{amsmath}    
\usepackage{amssymb}    
\usepackage{amsthm}     
\usepackage{booktabs}   
\usepackage{graphicx}   
\usepackage{bm}         
\usepackage{braket}     
\usepackage{longtable}  
\usepackage[margin=2.5cm]{geometry}
\usepackage{hhline}
\usepackage[dvipsnames]{xcolor}
\usepackage{changepage}

\usepackage[T1]{fontenc}
\usepackage{lmodern} 

\usepackage[colorlinks=true, citecolor=Purple, linkcolor=OliveGreen, urlcolor=blue, pdfborder={0 0 1}]{hyperref}

\usepackage{listings}

\usepackage[
    backend=biber,      
    sorting=none,
    autocite=inline,    
]{biblatex}
\DeclareFieldFormat{url}{\newline{\footnotesize\mbox{\url{#1}}}}

\AtBeginBibliography{%
}

\title{Enumeration of Polyominoes up to Size N=59}
\author{Toshihiro Shirakawa} 
\date{\today} 

\begin{document}

\nocite{Golomb1954Checker, Golomb_Polyominoes_JP, Lunnon1975, Read_1978_BoundsPoly, Redelmeier_1981_CountingPoly, OliveiraSilva_2007_N28, PolyominoesWeb, Shirakawa_G4G10_Polyomino, Mason_2023_Size50, Jensen_polyominoes_ICCS2003, Barequet_2022_N70, OEIS}

\maketitle 

\begin{abstract}
This paper reports the results of numerical computations for determining the number of polyominoes of size n (n-ominoes).
We verify the existing counts for $n \leq 50$ and newly compute the total number of polyominoes up to $n \leq 59$,
extending the counting limit. This work shows that, in addition to optimizing the search algorithm for the polyomino counting problem,
multi-threading dramatically improves computational efficiency.
\end{abstract}

\section{Introduction}
A \textbf{polyomino} is a connected figure formed by connecting \textbf{square cells} edge to edge.
A polyomino composed of $n$ cells is called an \textbf{$n$-omino}, and the total number of distinct $n$-ominoes is denoted by $\mathrm{Free}(n)$.
The problem of enumerating $n$-ominoes is one of the classic and \textbf{unsolved problems} in combinatorics.
While the exact value of $\mathrm{Free}(n)$ is known for small $n$, a \textbf{closed-form expression} for general $n$ has not yet been discovered.
The purpose of this paper is to report the values of $\mathrm{Free}(n)$ for $n$ that have not been previously computed, using the latest algorithms.

\section{Definitions and Methods of Enumeration}

\subsection{Definition of Polyominoes}
The enumeration in this paper focuses on the number of \textbf{Free Polyominoes} (OEIS:A000105), where polyominoes that are congruent by rotation and reflection are considered identical.
Hereafter, unless otherwise noted, the term ``polyomino'' refers to a Free Polyomino.
We also report the number of \textbf{One-sided Polyominoes} (OEIS:A000988), where polyominoes are considered identical only if they are congruent by rotation, and mirror images by reflection are distinguished.
Furthermore, the case where all rotations and reflections are distinguished, and only translation leads to identification, is referred to as \textbf{Fixed Polyominoes}.

\subsection{History of Polyomino Enumeration}
Since the concept of polyominoes was devised by Solomon Golomb in 1954 \cite{Golomb1954Checker}, their numbers have been sought by various mathematicians and puzzle enthusiasts.
While early search results are difficult to track as consolidated information, one example of the remaining results is as follows.\\
\begin{minipage}{\linewidth} %
\begin{tabular}{ll}
$n\leq 16$ & William Lunnon(1975)\cite{Lunnon1975}\\
$17\leq n \leq18$ & David R. Read(1978)\cite{Read_1978_BoundsPoly}\\
$19\leq n\leq 24$ & D. Hugh Redelmeier(1981)\cite{Redelmeier_1981_CountingPoly}\\
$n=25$ & Conway(1995)\cite{Jensen_polyominoes_ICCS2003}\\
$26\leq n\leq 28$ & Iwan Jensen(1998-2002)\cite{Jensen_polyominoes_ICCS2003}, Tomás Oliveira e Silva(2002)\cite{OliveiraSilva_2007_N28}\\
$29\leq n\leq 45,n=47$ & Toshihiro Shirakawa(2009)\cite{PolyominoesWeb}\cite{Shirakawa_G4G10_Polyomino}\\
$n=46,48\leq n\leq 50$ & John Mason(2021-2023)\cite{Mason_2023_Size50}\\
\end{tabular}
\end{minipage}
Furthermore, the count of Fixed Polyominoes was significantly updated by the enumeration for $n\leq 56$ using the \textbf{transfer-matrix algorithm} by Jensen \cite{Jensen_polyominoes_ICCS2003},
and further updated for $57\leq n\leq 70$ by Barequet and Ben-Shachar \cite{Barequet_2022_N70}, who improved the transfer-matrix algorithm.
\subsection{Overview of the Computation}
Similar to the calculation performed in 2009 \cite{Shirakawa_G4G10_Polyomino}, we determine the total number of polyominoes by using Burnside's Lemma to efficiently handle symmetry classes.
The total count can be determined by enumerating four types of polyominoes: \textbf{Fixed Polyominoes}, \textbf{Mirror-symmetric Polyominoes}, \textbf{Point-symmetric Polyominoes}, and \textbf{90-degree Rotationally Symmetric Polyominoes}.
We utilized the known calculation results for Fixed Polyominoes and performed new calculations for the remaining types.
Note that when applying Burnside's Lemma, we count all polyominoes that possess at least that specific symmetry, regardless of whether they possess higher symmetry.
The names for each symmetry are based on Mason's nomenclature \cite{Mason_2023_Size50}. While we count those with higher symmetry, we will use the same symbols for notation.
Compared to the method of listing polyominoes with exactly that symmetry, this approach offers two advantages: first, it avoids the need to determine whether a higher symmetry exists; and second, it eliminates the need to separately search for bi-axially symmetric shapes.
When categorized by the position of the axis or center of symmetry, there are a total of eight cases:
\begin{itemize}
    \item \textbf{Mirror symmetry} (3 types): vertical through cell boundaries (M90), vertical through cell centers (M90V), and diagonal (M45).
    \item \textbf{Point symmetry} (3 types): center on a cell center (R180C), center on an edge center (R180M), and center on a vertex (R180V).
    \item \textbf{90-degree rotational symmetry} (2 types): center on a cell center (R90C), and center on a vertex (R90V).
\end{itemize}
Since the number of M90V polyominoes is equal to $\mathrm{Fixed}(n/2)$, only the remaining seven types are subject to exploration.
The formulas for determining the number of Free Polyominoes ($\mathrm{Free}(n)$) and One-sided Polyominoes ($\text{One-sided}(n)$) from the calculated counts for each symmetry are as follows. The term $\mathrm{Fixed}(n/2)$ is set to 0 when $n$ is odd.
\begin{flalign}
\label{eq:free}
\mathrm{Free}(n) = \frac{1}{8}\cdot ( & \mathrm{Fixed}(n)+2\cdot \mathrm{Fixed}(n/2)+2\cdot \mathrm{M90}(n)+2\cdot \mathrm{M45}(n)\nonumber \\
& +\mathrm{R180C}(n)+2\cdot \mathrm{R180M}(n)+\mathrm{R180V}(n)+2\cdot \mathrm{R90C}(n)+2\cdot \mathrm{R90V}(n))&
\end{flalign}
\begin{flalign}
\label{eq:one-sided}
\text{One-sided}(n)=\frac{1}{4}\cdot ( & \mathrm{Fixed}(n)+\mathrm{R180C}(n)+2\cdot \mathrm{R180M}(n)+\mathrm{R180V}(n)\nonumber \\
& +2\cdot \mathrm{R90C}(n)+2\cdot \mathrm{R90V}(n))&
\end{flalign}
\begin{figure}[h]
  \centering
  \includegraphics[width=10cm]{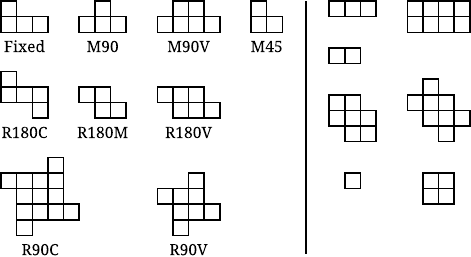}
  \caption{Symmetries used (Left) and unused (Right) in this calculation}
  \label{fig:Sym}
\end{figure}

\subsection{Detailed Algorithm}
\subsubsection{Mirror-symmetric Polyominoes}
\begin{figure}[t]
    \centering 
    \begin{minipage}[t]{0.35\linewidth}
        \centering
        \includegraphics[width=\linewidth]{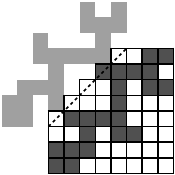}
        \caption{Example of M45. The dark gray area is the search domain.}
        \label{fig:Diag}
    \end{minipage}
    \hspace{2em} 
    \begin{minipage}[t]{0.35\linewidth}
        \centering
        \includegraphics[width=\linewidth]{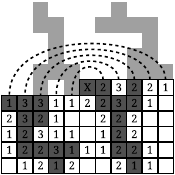}
        \caption{Example of R180C. Numbers indicate the values of the neighborhood counter.}
        \label{fig:Rot}
    \end{minipage}
\end{figure}

In this research, we adopted the \textbf{transfer-matrix algorithm} for the enumeration of mirror-symmetric polyominoes. For state encoding, we utilized the encoding based on Motzkin paths detailed in Jensen's doctoral thesis \cite{Jensen_Thesis_TM}.
The transfer-matrix method is applied to the region on one side of the axis of symmetry.
For M90 symmetry, we performed DP (Dynamic Programming) using the transfer-matrix method for both the vertical and horizontal orientations.
Since cells on the axis of symmetry are counted with a weight of 1, and all other cells with a weight of 2, a search conducted with a short side of the bounding box of length $W$ yields the correct value for $n \leq 3 \cdot W$. \\
For M45 symmetry, we applied the Transfer-Matrix Algorithm to the region defined by the bounding box on one side of the axis of symmetry.
The calculation was performed for each possible position where the axis intersects the bounding box. The required width of the bounding box varies depending on the number of cells the axis passes through, but the correct values were determined for up to $W=n/3$.
Figure \ref{fig:Diag} shows a diagonally mirror-symmetric polyomino, where the axis of symmetry passes through 5 cells of an $8 \times 8$ bounding box.
The total size of the symmetric polyomino (n) is enumerated by assigning a weight of 1 to the cells lying on the axis of symmetry and a weight of 2 to all other cells.
While the strict condition is that at least one cell must exist on all four sides and the diagonal, we only specified the maximum bounding box in this implementation. Therefore, a polyomino was considered valid if at least one cell existed on the top side, the left side, and the diagonal.
Since a simple implementation was sufficiently fast for $n \leq 60$, no active \textbf{Pruning} and \textbf{Parallelization} was performed.

\subsubsection{Point-symmetric Polyominoes}
For R180C, R180M, and R180V, we used the counting method proposed by Redelmeier \cite{Redelmeier_1981_CountingPoly}, which first calculates the rings and then enumerates the polyominoes belonging to those rings.
For the degenerate rings $1 \times 1$, $1 \times 2$, and $2 \times 2$, we adopted the search space partitioning method used by Mason \cite{Mason_2023_Size50}.
The implementation of this search space partitioning required the most time, and achieving a reduction in the \textbf{constant factor} of the computational complexity is crucial here. The source code for this optimization is provided in \textbf{Code Listing \ref{lst:optimizecode}}.
This code operates as a search for R180C $1 \times 1$, but it does not include search space reduction due to symmetry and it is not implemented with multi-threading.
The detailed search algorithm is as follows.
First, the starting cell (cell 'X' in Figure \ref{fig:Rot}) is marked as visited, and the unvisited cells among its four neighbors are added to the frontier.
A table for the Neighborhood Counter is prepared, and when a cell is added to the polyomino, the counter for its four neighbors is checked.
If the counter is 0, the cell is added to the frontier; regardless of whether it was 0, the Neighborhood Counter is incremented by +1.
Conversely, when a cell is removed from the polyomino during backtracking, the Neighborhood Counter for all four neighbors is decremented by -1.
This mechanism ensures that a neighboring cell is added to the frontier only the first time it is adjoined by a cell added to the polyomino.
Furthermore, we explore only one side of the axis of symmetry, and cells that extend beyond the boundary wrap around to the opposite side.
This is implemented within the function that enumerates the four neighbors, making cells connected by dashed lines behave as if their top edges are adjacent to each other.
The backtracking part was optimized by being implemented with a loop instead of recursion, and further accelerated by omitting and aggregating the final two stages.
For parallelization, we used a technique where, upon reaching a specific depth, the number of times that depth has been reached is counted, and the subsequent search is performed only if the count modulo the number of threads matches the thread ID.
This approach avoids the need for synchronization between threads, at the cost of requiring duplicated exploration of the shallow search space.

\subsubsection{90-degree Rotationally Symmetric Polyominoes}
Since the results for R90C (OEIS:A348849) and R90V (A348848) are already known, we performed a re-verification. Since the computation time is short, no special optimization was implemented.
\newpage
\vspace*{-1.5cm}
\section{Individual Calculation Results}
\begin{center}
\footnotesize
\renewcommand{\arraystretch}{0.9}
{\fontfamily{ntxtlf}\selectfont
\begin{tabular}{|r|r|r|r|}
\hhline{|-|-|-|-|}
n & M90 (OEIS:A346799) & M45 (OEIS:A346800) & R180C \\ \hhline{|-|-|-|-|}
1 & 1 & 1 & 1 \\ \hhline{|-|-|-|-|}
2 & 1 & 0 & 0 \\ \hhline{|-|-|-|-|}
3 & 2 & 2 & 2 \\ \hhline{|-|-|-|-|}
4 & 3 & 1 & 0 \\ \hhline{|-|-|-|-|}
5 & 7 & 5 & 7 \\ \hhline{|-|-|-|-|}
6 & 10 & 4 & 0 \\ \hhline{|-|-|-|-|}
7 & 24 & 16 & 24 \\ \hhline{|-|-|-|-|}
8 & 36 & 13 & 1 \\ \hhline{|-|-|-|-|}
9 & 86 & 54 & 86 \\ \hhline{|-|-|-|-|}
10 & 133 & 46 & 6 \\ \hhline{|-|-|-|-|}
11 & 314 & 186 & 316 \\ \hhline{|-|-|-|-|}
12 & 499 & 167 & 27 \\ \hhline{|-|-|-|-|}
13 & 1164 & 660 & 1178 \\ \hhline{|-|-|-|-|}
14 & 1888 & 612 & 114 \\ \hhline{|-|-|-|-|}
15 & 4366 & 2384 & 4442 \\ \hhline{|-|-|-|-|}
16 & 7192 & 2267 & 472 \\ \hhline{|-|-|-|-|}
17 & 16522 & 8726 & 16896 \\ \hhline{|-|-|-|-|}
18 & 27548 & 8464 & 1926 \\ \hhline{|-|-|-|-|}
19 & 62954 & 32278 & 64716 \\ \hhline{|-|-|-|-|}
20 & 106004 & 31822 & 7789 \\ \hhline{|-|-|-|-|}
21 & 241203 & 120419 & 249283 \\ \hhline{|-|-|-|-|}
22 & 409492 & 120338 & 31348 \\ \hhline{|-|-|-|-|}
23 & 928376 & 452420 & 964708 \\ \hhline{|-|-|-|-|}
24 & 1587151 & 457320 & 125802 \\ \hhline{|-|-|-|-|}
25 & 3586999 & 1709845 & 3748031 \\ \hhline{|-|-|-|-|}
26 & 6169400 & 1745438 & 503948 \\ \hhline{|-|-|-|-|}
27 & 13904736 & 6494848 & 14610484 \\ \hhline{|-|-|-|-|}
28 & 24041597 & 6686929 & 2016677 \\ \hhline{|-|-|-|-|}
29 & 54053950 & 24779026 & 57118440 \\ \hhline{|-|-|-|-|}
30 & 93896826 & 25703792 & 8065830 \\ \hhline{|-|-|-|-|}
31 & 210654990 & 94899470 & 223859532 \\ \hhline{|-|-|-|-|}
32 & 367450477 & 99096382 & 32251819 \\ \hhline{|-|-|-|-|}
33 & 822754494 & 364680344 & 879285686 \\ \hhline{|-|-|-|-|}
34 & 1440514144 & 383067646 & 128955260 \\ \hhline{|-|-|-|-|}
35 & 3219725534 & 1405619344 & 3460424846 \\ \hhline{|-|-|-|-|}
36 & 5656283859 & 1484352159 & 515657653 \\ \hhline{|-|-|-|-|}
37 & 12622055937 & 5432421429 & 13642112667 \\ \hhline{|-|-|-|-|}
38 & 22242057564 & 5764277096 & 2062335114 \\ \hhline{|-|-|-|-|}
39 & 49559836758 & 21046198560 & 53865266960 \\ \hhline{|-|-|-|-|}
40 & 87577573856 & 22429257682 & 8250061654 \\ \hhline{|-|-|-|-|}
41 & 194874338805 & 81716371069 & 212982833863 \\ \hhline{|-|-|-|-|}
42 & 345252222481 & 87432657722 & 33011955188 \\ \hhline{|-|-|-|-|}
43 & 767272498486 & 317917129256 & 843202561450 \\ \hhline{|-|-|-|-|}
44 & 1362583793862 & 341394729018 & 132132934138 \\ \hhline{|-|-|-|-|}
45 & 3024594142570 & 1239120776640 & 3342114787290 \\ \hhline{|-|-|-|-|}
46 & 5383141106732 & 1335080732960 & 529032164150 \\ \hhline{|-|-|-|-|}
47 & 11936178797526 & 4837744188806 & 13260804930438 \\ \hhline{|-|-|-|-|}
48 & 21287323864258 & 5228480834780 & 2118783780213 \\ \hhline{|-|-|-|-|}
49 & 47152797568454 & 18916648972476 & 52667141592300 \\ \hhline{|-|-|-|-|}
50 & 84253847090995 & 20502932765536 & 8488385208118 \\ \hhline{|-|-|-|-|}
51 & 186449386168092 & 74074144397536 & 209361911445176 \\ \hhline{|-|-|-|-|}
52 & 333745580076968 & 80498698956092 & 34017019523282 \\ \hhline{|-|-|-|-|}
53 & 737898749100214 & 290445682120726 & 832940890785598 \\ \hhline{|-|-|-|-|}
54 & 1323045715291324 & 316415797789354 & 136363468404368 \\ \hhline{|-|-|-|-|}
55 & 2922731717041198 & 1140246035002130 & 3316370384325374 \\ \hhline{|-|-|-|-|}
56 & 5248636433672399 & 1245063333064256 & 546798865890672 \\ \hhline{|-|-|-|-|}
57 & 11585483714367491 & 4481581094216159 & 13213587492371313 \\ \hhline{|-|-|-|-|}
58 & 20835850857774290 & 4904100642586164 & 2193217495584004 \\ \hhline{|-|-|-|-|}
59 & 45956970811133468 & 17633164855044618 & 52682486874469758 \\ \hhline{|-|-|-|-|}
60 & 82765943504776434 & 19334605766172446 &  \\ \hhline{|-|-|-|-|}
61 & & 69449002269748050 &  \\ \hhline{|-|-|-|-|}
62 & & 76294817933002026 &  \\ \hhline{|-|-|-|-|}
63 & & 273786260678453878 &  \\ \hhline{|-|-|-|-|}
64 & & 301311979812960930 &  \\ \hhline{|-|-|-|-|}
65 & & 1080295987815098059 &  \\ \hhline{|-|-|-|-|}
66 & & 1190910519195326830 &  \\ \hhline{|-|-|-|-|}
\end{tabular}
}
\end{center}

\newpage
\begin{center}
\footnotesize
\renewcommand{\arraystretch}{1}
{\fontfamily{ntxtlf}\selectfont
\begin{tabular}{|r|r|r|r|}
\hhline{|-|-|-|-|}
n & R180M $1\times 2$ core & R180M other rings & R180M total \\ \hhline{|-|-|-|-|}
2 & 1 & 0 & 1 \\ \hhline{|-|-|-|-|}
4 & 3 & 0 & 3 \\ \hhline{|-|-|-|-|}
6 & 10 & 0 & 10 \\ \hhline{|-|-|-|-|}
8 & 35 & 0 & 35 \\ \hhline{|-|-|-|-|}
10 & 126 & 1 & 127 \\ \hhline{|-|-|-|-|}
12 & 465 & 7 & 472 \\ \hhline{|-|-|-|-|}
14 & 1742 & 36 & 1778 \\ \hhline{|-|-|-|-|}
16 & 6594 & 166 & 6760 \\ \hhline{|-|-|-|-|}
18 & 25165 & 731 & 25896 \\ \hhline{|-|-|-|-|}
20 & 96651 & 3132 & 99783 \\ \hhline{|-|-|-|-|}
22 & 373138 & 13182 & 386320 \\ \hhline{|-|-|-|-|}
24 & 1446826 & 54835 & 1501661 \\ \hhline{|-|-|-|-|}
26 & 5630471 & 226277 & 5856748 \\ \hhline{|-|-|-|-|}
28 & 21979840 & 928390 & 22908230 \\ \hhline{|-|-|-|-|}
30 & 86034888 & 3793184 & 89828072 \\ \hhline{|-|-|-|-|}
32 & 337556423 & 15449979 & 353006402 \\ \hhline{|-|-|-|-|}
34 & 1327143622 & 62781732 & 1389925354 \\ \hhline{|-|-|-|-|}
36 & 5227430169 & 254659893 & 5482090062 \\ \hhline{|-|-|-|-|}
38 & 20623974712 & 1031544202 & 21655518914 \\ \hhline{|-|-|-|-|}
40 & 81488847870 & 4173969124 & 85662816994 \\ \hhline{|-|-|-|-|}
42 & 322406638993 & 16875126628 & 339281765621 \\ \hhline{|-|-|-|-|}
44 & 1277134313760 & 68180532058 & 1345314845818 \\ \hhline{|-|-|-|-|}
46 & 5064653939573 & 275328459109 & 5339982398682 \\ \hhline{|-|-|-|-|}
48 & 20104959952132 & 1111393767171 & 21216353719303 \\ \hhline{|-|-|-|-|}
50 & 79884293217360 & 4484868477645 & 84369161695005 \\ \hhline{|-|-|-|-|}
52 & 317682265951332 & 18093701200360 & 335775967151692 \\ \hhline{|-|-|-|-|}
54 & 1264358079284029 & 72983734261285 & 1337341813545314 \\ \hhline{|-|-|-|-|}
56 & 5035793376023652 & 294351310722927 & 5330144686746579 \\ \hhline{|-|-|-|-|}
58 & 20070775724211511 & 1187035156398203 & 21257810880609714 \\ \hhline{|-|-|-|-|}
\end{tabular}
}
\end{center}

\begin{center}
\footnotesize
\renewcommand{\arraystretch}{1}
{\fontfamily{ntxtlf}\selectfont
\begin{tabular}{|r|r|r|r|}
\hhline{|-|-|-|-|}
n & R180V $2\times 2$ core & R180V other rings & R180V total \\ \hhline{|-|-|-|-|}
2 & 0 & 0 & 0 \\ \hhline{|-|-|-|-|}
4 & 1 & 0 & 1 \\ \hhline{|-|-|-|-|}
6 & 4 & 0 & 4 \\ \hhline{|-|-|-|-|}
8 & 14 & 0 & 14 \\ \hhline{|-|-|-|-|}
10 & 50 & 0 & 50 \\ \hhline{|-|-|-|-|}
12 & 183 & 3 & 186 \\ \hhline{|-|-|-|-|}
14 & 682 & 22 & 704 \\ \hhline{|-|-|-|-|}
16 & 2569 & 116 & 2685 \\ \hhline{|-|-|-|-|}
18 & 9760 & 552 & 10312 \\ \hhline{|-|-|-|-|}
20 & 37347 & 2494 & 39841 \\ \hhline{|-|-|-|-|}
22 & 143744 & 10898 & 154642 \\ \hhline{|-|-|-|-|}
24 & 555876 & 46592 & 602468 \\ \hhline{|-|-|-|-|}
26 & 2158290 & 196302 & 2354592 \\ \hhline{|-|-|-|-|}
28 & 8408823 & 818589 & 9227412 \\ \hhline{|-|-|-|-|}
30 & 32858248 & 3388052 & 36246300 \\ \hhline{|-|-|-|-|}
32 & 128726908 & 13945033 & 142671941 \\ \hhline{|-|-|-|-|}
34 & 505444284 & 57156614 & 562600898 \\ \hhline{|-|-|-|-|}
36 & 1988576779 & 233515451 & 2222092230 \\ \hhline{|-|-|-|-|}
38 & 7837564226 & 951653244 & 8789217470 \\ \hhline{|-|-|-|-|}
40 & 30939196525 & 3870694267 & 34809890792 \\ \hhline{|-|-|-|-|}
42 & 122308756306 & 15718933882 & 138027690188 \\ \hhline{|-|-|-|-|}
44 & 484136456962 & 63755513030 & 547891969992 \\ \hhline{|-|-|-|-|}
46 & 1918619453466 & 258332542690 & 2176951996156 \\ \hhline{|-|-|-|-|}
48 & 7611608307942 & 1045901961521 & 8657510269463 \\ \hhline{|-|-|-|-|}
50 & 30226772614264 & 4231749935548 & 34458522549812 \\ \hhline{|-|-|-|-|}
52 & 120143630598288 & 17112745629450 & 137256376227738 \\ \hhline{|-|-|-|-|}
54 & 477939741578208 & 69172483040912 & 547112224619120 \\ \hhline{|-|-|-|-|}
56 & 1902751625571409 & 279509304959259 & 2182260930530668 \\ \hhline{|-|-|-|-|}
58 & 7580601028004440 & 1129112738593514 & 8709713766597954 \\ \hhline{|-|-|-|-|}
\end{tabular}
}
\end{center}

\newpage
\section{Enumeration Results}
The values for $n\leq 50$ in normal font were known prior to this study, and the values for $51\leq n\leq 59$ in \textbf{bold} are the newly computed numbers.
The section titled ``Individual Calculation Results'' contains the enumeration results for each symmetry, while ``Enumeration Results'' presents the final number of polyominoes, which is the main objective of this paper.
\begin{center}
\footnotesize
\renewcommand{\arraystretch}{0.9}
{\fontfamily{ntxtlf}\selectfont
\begin{tabular}{|r|r|r|}
\hhline{|-|-|-|} 
$n$ & Polyominoes (OEIS:A000105) & One-sided Polyominoes (OEIS:A000988) \\ \hhline{|-|-|-|}
1 & 1 & 1 \\ \hhline{|-|-|-|}
2 & 1 & 1 \\ \hhline{|-|-|-|}
3 & 2 & 2 \\ \hhline{|-|-|-|}
4 & 5 & 7 \\ \hhline{|-|-|-|}
5 & 12 & 18 \\ \hhline{|-|-|-|}
6 & 35 & 60 \\ \hhline{|-|-|-|}
7 & 108 & 196 \\ \hhline{|-|-|-|}
8 & 369 & 704 \\ \hhline{|-|-|-|}
9 & 1285 & 2500 \\ \hhline{|-|-|-|}
10 & 4655 & 9189 \\ \hhline{|-|-|-|}
11 & 17073 & 33896 \\ \hhline{|-|-|-|}
12 & 63600 & 126759 \\ \hhline{|-|-|-|}
13 & 238591 & 476270 \\ \hhline{|-|-|-|}
14 & 901971 & 1802312 \\ \hhline{|-|-|-|}
15 & 3426576 & 6849777 \\ \hhline{|-|-|-|}
16 & 13079255 & 26152418 \\ \hhline{|-|-|-|}
17 & 50107909 & 100203194 \\ \hhline{|-|-|-|}
18 & 192622052 & 385221143 \\ \hhline{|-|-|-|}
19 & 742624232 & 1485200848 \\ \hhline{|-|-|-|}
20 & 2870671950 & 5741256764 \\ \hhline{|-|-|-|}
21 & 11123060678 & 22245940545 \\ \hhline{|-|-|-|}
22 & 43191857688 & 86383382827 \\ \hhline{|-|-|-|}
23 & 168047007728 & 336093325058 \\ \hhline{|-|-|-|}
24 & 654999700403 & 1309998125640 \\ \hhline{|-|-|-|}
25 & 2557227044764 & 5114451441106 \\ \hhline{|-|-|-|}
26 & 9999088822075 & 19998172734786 \\ \hhline{|-|-|-|}
27 & 39153010938487 & 78306011677182 \\ \hhline{|-|-|-|}
28 & 153511100594603 & 307022182222506 \\ \hhline{|-|-|-|}
29 & 602621953061978 & 1205243866707468 \\ \hhline{|-|-|-|}
30 & 2368347037571252 & 4736694001644862 \\ \hhline{|-|-|-|}
31 & 9317706529987950 & 18635412907198670 \\ \hhline{|-|-|-|}
32 & 36695016991712879 & 73390033697855860 \\ \hhline{|-|-|-|}
33 & 144648268175306702 & 289296535756895985 \\ \hhline{|-|-|-|}
34 & 570694242129491412 & 1141388483146794007 \\ \hhline{|-|-|-|}
35 & 2253491528465905342 & 4506983054619138245 \\ \hhline{|-|-|-|}
36 & 8905339105809603405 & 17810678207278478530 \\ \hhline{|-|-|-|}
37 & 35218318816847951974 & 70436637624668665265 \\ \hhline{|-|-|-|}
38 & 139377733711832678648 & 278755467406691820628 \\ \hhline{|-|-|-|}
39 & 551961891896743223274 & 1103923783758183428889 \\ \hhline{|-|-|-|}
40 & 2187263896664830239467 & 4374527793263174673335 \\ \hhline{|-|-|-|}
41 & 8672737591212363420225 & 17345475182286431485513 \\ \hhline{|-|-|-|}
42 & 34408176607279501779592 & 68816353214298169362691 \\ \hhline{|-|-|-|}
43 & 136585913609703198598627 & 273171827218863802383383 \\ \hhline{|-|-|-|}
44 & 542473001706357882732070 & 1084946003411691009916361 \\ \hhline{|-|-|-|}
45 & 2155600091107324229254415 & 4311200182212516601049225 \\ \hhline{|-|-|-|}
46 & 8569720333296834568434605 & 17139440666589637839781602 \\ \hhline{|-|-|-|}
47 & 34085105553123831158180217 & 68170211106239275354867268 \\ \hhline{|-|-|-|}
48 & 135629410647775553284438364 & 271258821295535228672142075 \\ \hhline{|-|-|-|}
49 & 539916438668093786698843965 & 1079832877336154538674417465 \\ \hhline{|-|-|-|}
50 & 2150182610161041739167164220 & 4300365220322020871043392169 \\ \hhline{|-|-|-|}
51 & \bfseries{8566301646855786503391647670} & \bfseries{17132603293711442745018012526} \\ \hhline{|-|-|-|}
52 & \bfseries{34140832502065121824070942627} & \bfseries{68281665004129996529664185170} \\ \hhline{|-|-|-|}
53 & \bfseries{136116710694919148806231527457} & \bfseries{272233421389837783440247444444} \\ \hhline{|-|-|-|}
54 & \bfseries{542874683702992552434145546370} & \bfseries{1085749367405984128525518503279} \\ \hhline{|-|-|-|}
55 & \bfseries{2165873377717239478407665345440} & \bfseries{4331746755434476925326454669216} \\ \hhline{|-|-|-|}
56 & \bfseries{8643839320316615721187562237534} & \bfseries{17287678640633227581480905193254} \\ \hhline{|-|-|-|}
57 & \bfseries{34507662858721137395351058631991} & \bfseries{69015325717442266757169712972157} \\ \hhline{|-|-|-|}
58 & \bfseries{137801831798576866546096116428465} & \bfseries{275603663597153717811728777821645} \\ \hhline{|-|-|-|}
59 & \bfseries{550453451340976338599795923539996} & \bfseries{1100906902681952645404524013990949} \\ \hhline{|-|-|-|}
\end{tabular}
}
\end{center}

\newpage
\section{Computation Time}
This section details the performance of the computer used and the actual time taken for the computations.

\subsection{Computational Efficiency and Contribution}
This research demonstrates the high-speed enumeration of mirror-symmetric polyominoes based on the transfer-matrix algorithm established by Jensen.
For point-symmetric polyominoes with an even number of cells, the use of Redelmeier's method of first calculating the rings significantly improved the computation time compared to the $n\leq 45$ calculation from 2009.
For the tasks requiring sequential search, specifically R180C $1\times 1$, R180M $1\times 2$, and R180V $2\times 2$, no order-level improvement was made, but implementation-level optimizations such as the reduction of the search space via symmetry and multi-threading were applied.
As a result, a speedup of approximately 5000 times was achieved even for the most time-consuming sequential search tasks. However, these tasks remain the most time-intensive due to the still large number of objects to be searched.

\subsection{Execution Environment and Results}
The following computer was used for the calculations.
\begin{table}[ht]
    \centering
    \begin{tabular}{lll} \hline
    Motherboard & Supermicro H11DSi rev2 \\
    CPU & AMD EPYC 7452 (2.35\,GHz) $\times2$ & 64C128T \\
    Memory & DDR4-2666 ECC 64\,GB $\times8$ & Total 512 GB \\ \hline
    \end{tabular}
    \caption{Computer Specifications}
    \label{table:PC}
\end{table}

Since a simple implementation was sufficiently fast for $n\leq 60$, no active Pruning and Parallelization was performed for M90 and M45.
Consequently, the memory limit was reached at $n\leq 60$ for M90 and $n\leq 66$ for M45.
However, it is expected that the maximum solvable $n$ could be significantly increased by applying pruning based on connectivity and parallelization.
There are also slight differences in the level of optimization for the other tasks, suggesting that implementation choices could result in several tens of percent difference in performance.

\begin{table}[ht]
    \centering
    \begin{tabular}{|l|r|r|r|r|l|}
    \hhline{|-|-|-|-|-|-|}
    \textbf{Task} & \textbf{Hours} & $\mathbf{N_{max}}$ & \textbf{Count(E12)} & \textbf{Threads} & \textbf{Remarks} \\ \hhline{|-|-|-|-|-|-|}
    M90 & 10 & 60 & 82765 & 1 & \\ \hhline{|-|-|-|-|-|-|}
    M45 & 6 & 60 & 19334 & 1 & 30h for $n\leq 66$ \\ \hhline{|-|-|-|-|-|-|}
    R180C $1\times 1$ & 270 & 59 & 52682 & 128 & \\ \hhline{|-|-|-|-|-|-|}
    R180C other & 19 & 58 & 2193 & 128 & \\ \hhline{|-|-|-|-|-|-|}
    R180M $1\times 2$ & 138 & 58 & 20070 & 128 & \\ \hhline{|-|-|-|-|-|-|}
    R180M other & 8 & 58 & 1187 & 128 &  \\ \hhline{|-|-|-|-|-|-|}
    R180V $2\times 2$ & 30 & 58 & 7580 & 128 & \\ \hhline{|-|-|-|-|-|-|}
    R180V other & 5 & 58 & 1129 & 128 & \\ \hhline{|-|-|-|-|-|-|}
    R90C & $<$1 & 61 & $<$1 & 1 & \\ \hhline{|-|-|-|-|-|-|}
    R90V & $<$1 & 60 &$<$1 & 1 & \\ \hhline{|-|-|-|-|-|-|}
    Total & 487 & 59 & & & $\approx 284$ hours for $n\leq 58$ \\ \hhline{|-|-|-|-|-|-|}
    \end{tabular}
    \caption{Computation Time}
    \label{tab:time}
\end{table}


\section{Conclusion}
The computed values are in complete agreement with previously reported values ($n \leq 50$).
This agreement not only confirms the validity of the algorithm implementation used in this paper but also verifies the previous results.
Therefore, the newly computed values for $n \leq 59$ reported in this study are highly reliable.

\subsection{Future Work}
The majority of the current computation time is consumed by point-symmetric polyominoes that have a core at the center.
It is suggested that the transfer-matrix method could be applicable to point-symmetric shapes with a central core,
and its application is expected to lead to further reduction in computation time and new record-breaking enumerations.

\newpage
\section{Code Listings} 
\renewcommand{\lstlistingname}{Code} 
\vspace{1ex}
\begin{lstlisting}[
language=C++,
caption={Optimized search code in C++},
label={lst:optimizecode}
]
#include <cstdio>
#include <array>
#include <ctime>
constexpr int MAXN = 61; // R180C
// constexpr int MAXN = 62; // R180M, R180V
long long calc(int n){
    unsigned char field[MAXN * MAXN / 2] = {};
    auto Neighborhood = [&](int pos) { // Enumerate neighboring cells
        std::array<int, 4> nbhs;
        nbhs[0] = pos + 1;
        nbhs[1] = pos + MAXN;
        nbhs[2] = pos - 1;
        nbhs[3] = pos - MAXN;
        nbhs[3] *= (nbhs[3] >> 31) | 1; // R180C
        // nbhs[3] ^= nbhs[3] >> 31; // R180M, R180V
        return nbhs;};
    n /= 2;
    std::array<int, MAXN * 2> frontier;
    std::array<int, MAXN / 2> put;
    int putSize = 0;
    field[0] = field[1] = field[MAXN] = 1;
    frontier[0] = 1;
    frontier[1] = MAXN;
    int frontierSize = 2;
    int frontierIndex = 0;
    long long sol = 0;
    while (1) { // Search loop
        if (frontierSize == frontierIndex) { // Backtrack
            if (putSize == 0) { // Search end
                break;}
            putSize--;
            int back = put[putSize];
            frontierIndex = back + 1;
            auto nbhs = Neighborhood(frontier[back]);
            for (const auto nbh : nbhs) {
                field[nbh]--;
                frontierSize -= !field[nbh];
        }} else {
            if (putSize == n - 2) { // Tally count for tail-end solutions
                auto d = frontierSize - frontierIndex;
                sol += d * (d - 1) / 2;
                while (frontierIndex < frontierSize) {
                    auto pos = frontier[frontierIndex];
                    auto nbs = Neighborhood(pos);
                    for (const auto nb : nbs) {
                        sol += !field[nb];}
                    frontierIndex++;
            }}else{ // Place cell
                auto pos = frontier[frontierIndex];
                auto nbhs = Neighborhood(pos);
                for (const auto nbh : nbhs) {
                    frontier[frontierSize] = nbh;
                    frontierSize += !field[nbh];
                    field[nbh]++;}
                put[putSize] = frontierIndex;
                putSize++;
                frontierIndex++;
    }}}
    return sol;
}
int main(){
    for (int i = 5; i < 60; i += 2){
        auto start = time(NULL);
        auto sol = calc(i);
        auto end = time(NULL);
        printf("%2d:%6ds, %lld\n", i, (int)(end - start), sol);}
    return 0;
}
\end{lstlisting}
\clearpage
\begin{adjustwidth}{-3mm}{-3mm}
\printbibliography
\end{adjustwidth}

\end{document}